\documentclass[11pt]{article}
\usepackage{amsfonts}
\usepackage{mathrsfs}
\usepackage{amsmath}
\usepackage{amssymb}
\usepackage{graphicx}
\graphicspath{{figure/}}
\usepackage{epic}
\usepackage{amsthm,latexsym}
\renewcommand{\paragraph}{\roman{paragraph}}
\newtheorem{theorem}{\scshape \mdseries  Theorem}[section]

\newtheorem{lemma}[theorem]{\scshape \mdseries  Lemma}
\newtheorem{coro}[theorem]{\scshape \mdseries  Corollary}
\newtheorem{conj}[theorem]{\scshape \mdseries  Conjecture}

\newtheorem{defi}[theorem]{\scshape \mdseries  Definition}

\begin{document}

\title{\sf On automorphism groups of power semigroups\\ over numerical semigroups  }
\author{ Dein Wong \thanks{Corresponding author, E-mail address:wongdein@163.com. Supported by  NSFC of China (No.12371025). School of Mathematics, China University of Mining and Technology, Xuzhou,  China. },\ \ \ \ \
 \ Songnian Xu \thanks{School of Mathematics, China University of Mining and Technology, Xuzhou,  China. },\ \ \ \ \ Chi Zhang \thanks{School of Mathematics, China University of Mining and Technology, Xuzhou,  China.},\ \ \ \
Jinxing Zhao\thanks{Corresponding author, E-mail address:. Supported by  NSFC of China (No.12161062).  School of Mathematics Sciences, Inner Mongolia University, Hohhot, China.}\ \ \ \ }
\date{}
\maketitle
\noindent {\bf Abstract:}\ \ A numerical semigroup  $S$ is  a cofinite subsemigroup   of $ \mathbb{N}$, where $\mathbb{N}$ is the additive monoid of non-negative integers. Denote by $\mathcal{P}_{\rm fin} (S)$ the semigroup consisting of all non-empty finite subsets of $S$ endowed with the operation of setwise addition defined by
$$X+Y=\{x+y:x\in X, y\in Y\}, \qquad\text{for all } X, Y \in \mathcal P_\text{fin}(S).$$ We call $\mathcal{P}_{\rm fin} (S)$ the finitary power semigroup of $S$.
 When $0\in S$ (and hence $S$ is a numerical monoid),      the family $\mathcal P_{\text{fin},0}(S)$ of all finite subsets of $S$ containing $0$ is  a submonoind of $\mathcal P_\text{fin}(S)$; we call   $\mathcal{P}_{{\rm fin}, 0}(S)$ the reduced finitary power monoid of $S$ with the singleton  $\{0\}$ as  zero-element.
  For a non-empty finite subset $X$ of $\mathbb{N}$, we denote by $ \min X$ and $\max X $ the minimum  and the maximum    in $X$.  Tringali and Yan have recently proved in [J.\ Combin.\ Theory Ser.\ A 209 (2025)]   that the only non-trivial automorphism of $\mathcal{P}_{{\rm fin},0}(\mathbb{N})$
is the involution $X \mapsto \max  X - X$.  By applying  Tringali-Yan's result, we in this article  determined the automorphism group of the finitary  power semigroup $\mathcal{P}_{\rm fin}(S)$ of an arbitrary numerical
semigroup $S$. More precisely,  if $S$ is the set of all integers larger than or equal to a fixed $k \in \mathbb N$, then the only non-trivial automorphism of $\mathcal{P}_{\rm fin}(S)$
is the involution $X \mapsto \max X - X+ \min X$; otherwise,  $\mathcal{P}_{\rm fin}(S)$ has only the identity automorphism.

\section{Introduction}
\quad \quad Let $H$ be a  semigroup. Denote by $\mathcal{P}_{\rm fin}(H)$ the  (finitary) power semigroup of
$H$, consisting of all  finite non-empty subsets of $H$ and endowed with the setwise binary operation $$(X, Y) \mapsto  XY := \{xy :\  x \in X, y \in  Y\}, \quad {\rm for\ all} \ X, Y \in \mathcal P_\text{fin}(H).\ $$
Moreover, if $H$ is a monoid with the  identity $1_H$,  we denote by $\mathcal{P}_{{\rm fin},1}(H)$ the set of finite subsets of $H$   containing    $1_H$;
 it is a submonoid of $\mathcal{P}_{\rm fin}(H)$ with identity $\{1_H\}$, henceforth called the reduced (finitary) power monoid of $H$.

Power semigroups  were first systematically studied by Shafer and Tamura \cite{sha1} in
the late 1960s. A central question of this field is the so-called ``Isomorphism Problem
For Power Semigroups": {\it Whether, for semigroups $S$ and $T$ in a certain class $\mathcal{O}$, an isomorphism
between  $\mathcal{P}_{\rm fin}(S)$ and  $\mathcal{P}_{\rm fin}(T)$ implies that $S$ and $T$ are isomorphic?}  Although this was answered in the
negative by Mogiljanskaja \cite{mog} for the class of all semigroups,  several other classes have been found for which the answer is positive.
Recently, power semigroups  were  investigated from multiple new perspectives   in a series of papers,  such as primality and atomicity \cite{agg};  arithmetic property \cite{ant};    algebraic properties  \cite{bie, dan};  factorization property     \cite{cos, cos1, fan, tri}.

More recently, some attention was concentrated on automorphism groups of power monoids or reduced power monoids of certain additive semigroups. Tringali and Yan \cite{tri2}   proved that:
 \begin{theorem}  {\rm (\cite{tri2}, Theorem 3.2)} The only   automorphisms of $\mathcal{P}_{\rm fin,0}(\mathbb{N})$ is either the identity or the
 involution $\sigma_0: X \mapsto \max X - X$, where $\max X$ is the maximum in $X$.\end{theorem}
 The involution  $\sigma_0$ is interesting in the following sense. Every automorphism $f$ of a monoid $H$ can be canonically
extended to an automorphism of $\mathcal{P}_{\rm fin}(H)$  or $\mathcal{P}_{{\rm fin}, 1}(H)$,
$$ X \mapsto f[X]:= \{f(x): x \in X\},\ \ \  {\rm for\ all}\   X\in  \mathcal{P}_{\rm fin}(H) \   {\rm or}\ X\in \mathcal{P}_{\rm fin,1}(H),$$ which is referred to  an inner automorphism  of $\mathcal{P}_{\rm fin}(H)$  or $\mathcal{P}_{\rm fin,1_H}(H)$. Thus the question: {\it Whether $\mathcal{P}_{\rm fin}(H)$  (resp., $\mathcal{P}_{{\rm fin},1}(H)$) has non-inner automorphisms?} is raised naturally. One can see that  the involution $\sigma_0: X \mapsto \max X - X$ for all $X\in  \mathcal{P}_{{\rm fin},0}(\mathbb{N})$ is not an inner automorphism  of
$\mathcal{P}_{{\rm fin},0}(\mathbb{N})$.
 Following up on Theorem 1.1, Tringali and Wen \cite{triwen} determined the automorphism group of the power monoid  $\mathcal{P}_{\rm fin}(\mathbb{Z})$,
 consisting of all finite subsets of the additive group of integers, and Rago \cite{rag}  gave a full description of the
automorphism group of $\mathcal{P}_{\rm fin,0}(G)$ for  a finite abelian group $G$.
Recall that a  numerical  semigroup (resp., a numerical monoid)  is a sub-semigroup (resp., submonoid)  $S$ of $(\mathbb{N }, +)$ such that $\mathbb{N}  \setminus S$ is a finite set. The present article is motivated by the following conjecture.
\begin{conj} {\rm (\cite{tri2}, Conjecture)}  The automorphism group of the reduced power monoid $\mathcal{P}_{\rm fin,0}(S)$ of a numerical
monoid $S$ properly contained in $\mathbb{N}$  must be  the identity.\end{conj}

  In this article, we do not address Tringali and Yan's conjecture, but consider an analogous problem: {\it How about the automorphism group of the   power monoid $\mathcal{P}_{\rm fin}(S)$ of a numerical
semigroup $S$ properly contained in $\mathbb{N}$?  whether or not it  contains only  the identity automorphism?}
In this article, we give this problem a definitely  answer.

\begin{theorem} If   $S$ is the set of all integers larger than or equal to a fixed $k \in \mathbb N$, then the only non-trivial automorphism of $\mathcal{P}_{\rm fin}(S)$
is the involution $\sigma: X \mapsto \max X -X+\min X$ for all $ X\in \mathcal{P}_{\rm fin}(S)$; otherwise, if $S$ is not a discrete interval, then $\mathcal{P}_{\rm fin}(S)$ has only the identity automorphism, where $\min X$  and $\max X$  are respectively the minimum   and the maximum    in $X$.
\end{theorem}

Two remarks should  be pointed out (as below):
\vskip 2mm $\bullet$ \ \ For the case when $S =[\![k,\infty):=\{ l\in \mathbb{N}:l\geq k\}  $ is a discrete interval, the involution $\sigma: X \mapsto \max X -X+\min X$ is not an inner automorphism of $\mathcal{P}_{\rm fin}(S)$, because the semigroup $S=[\![k,\infty)$ has only the identity automorphism  (see \cite{sas}, Corollary 1.4).

\vskip 2mm $\bullet$ \ \ Theorem 1.3 gives a characterization for the automorphism group of  $\mathcal{P}_{\rm fin}(S)$, which has, in principle, little to nothing to do with the
problem of determining the automorphism group   of the reduced finitary
power monoid $\mathcal{P}_{\rm fin,0}(S)$. Thus, Conjecture 1.2  still keeps open.
 \vskip 2mm

The technique for proving Theorem 1.3 is as follows:
\vskip 2mm
\noindent{\sf Step 1.}  Prove that $ \mathcal{P}_{\rm fin}([\![k,\infty))$ is stabilized by any automorphism $f$ of  $ \mathcal{P}_{\rm fin}(S)$, where   $k$ is the least member in $S$  such that $[\![k,\infty)\subseteq  S$, thus the problem of determining $f$ can be reduced to the restriction of $f$ to $ \mathcal{P}_{\rm fin}([\![k,\infty))$  (see Section 2).
\vskip 2mm
\noindent{\sf Step 2.}  Determine the automorphism group of $ \mathcal{P}_{\rm fin}([\![k,\infty))$  (see Section 3).
\vskip 2mm
\noindent{\sf Step 3.}
Based on the  form of an arbitrary automorphism of  $ \mathcal{P}_{\rm fin}([\![k,\infty))$,  complete the characterization   of any automorphism $f$ of  $ \mathcal{P}_{\rm fin}(S)$  (see Section 4).

\section{Reduce the problem to a special case}
\quad\quad We begin with some definitions and  some preliminaries that will be used for proving our main result. Write $\mathbb{N}$  for the set of non-negative integers. It forms a monoid under the  ordinary additive operation of integers and $0$ is its zero-element.  If $0\leq i\leq j \in \mathbb{N}$, we define  $[\![i ,j ]\!] $
to be the discrete interval $\{x \in\mathbb{N}: i \leq  x \leq  j\}$.  Given $l\in \mathbb{N}$  and $X\subseteq \mathbb{N}$, we set
$$\begin{array}{ccc}l+ X&=& \{l+x: x\in X\};\ \ \  \ \ \ \ \ \ \ \ \ \      \\
 l- X&=& \{l-x: x\in X\}. \ \ \ \ \ \ \ \ \ \ \ \ \
\end{array}$$
When $l$ is positive, by $lX$ we denote the $l$-fold sum of $X$, that is  $$lX= \{x_1 + \cdots + x_l : x_1, \ldots, x_l\in X\}.$$

For a numerical  semigroup $S$ of $\mathbb{N}$ with $S\not=\mathbb{N}$, there exists a unique positive integer $k$ such that $x\in S$ for all $k\leq x\in \mathbb{N}$ and $k-1\notin S$, which   is called the {\sf critical} element of $S$ and is written as $\theta_S$. The {\sf Frobenius number}  of $S$ is defined as $F(S)=: {\rm max}(\mathbb{N} \setminus S)$ (see \cite{bie}). Thus $\theta_S=F(S)+1$ if $S\not=\mathbb{N}$.
  We write $\alpha_S$ for the minimum in $S$. Then $0\leq \alpha_S\leq \theta_S$, $ \alpha_S=0$ implies that $S$ is a monoid and $\alpha_S=\theta_S$ implies that $S=[\![\alpha_S,\infty)$ is a discrete interval.
The  minimum   and the maximum    of a non-empty finite subset $X$ of $S$ are written as $\min X $ and $\max X$, respectively. For example, $S=\{0,3,5,6\}\cup [\![8,\infty]\!]$ is a   numerical  semigroup with $\alpha_S=0$ and $\theta_S=8$, for the subset $X=\{0,5,8,10\}$, we have $\min X=0$ and $\max X=10$.
 It is clear that \begin{eqnarray}\min (X+Y)=\min X+\min Y,\ \ \ \max(X+Y)=\max X+\max Y,\ \ \  {\rm for\ all}\  X, Y\in \mathcal{P}_{\rm fin}(S). \end{eqnarray}

 Hereafter, let $f$ be a given automorphism of $\mathcal{P}_{\rm fin}(S)$. The image of an element $X\in    \mathcal{P}_{\rm fin}(S)$ under $f$ is written as $X^f$.  Clearly, \begin{eqnarray} (X+Y)^f=X^f+Y^f,\ \ \ (l X)^f=l  X^f,\ \ \  {\rm for\ all}\  X \in \mathcal{P}_{\rm fin}(S),\ l\in \mathbb{N}\setminus\{0\}. \end{eqnarray}
To determine the form of the given automorphism $f$ of   $\mathcal{P}_{\rm fin}(S)$, a key technique is to find more enough fixed points of $f$ in  $\mathcal{P}_{\rm fin}(S)$.
It is easy to see that each singleton   is a   fixed point  of $f$.

\begin{lemma} Each  singleton $\{s\}$ of $\mathcal{P}_{\rm fin}(S)$, with $s\in S$, is    a   fixed point  of $f$.
\end{lemma}
\begin{proof}  It follows from  \cite{tri1}  (Theorem 1)
that  every automorphism of $\mathcal{P}_{\rm fin}(S)$ restricts to an
automorphism of $S$. Since the automorphism group of $S$ is trivial (see \cite{sas}, Corollary 1.4),
this immediately implies that $\{s\}$ is a fixed point of $f$ for each $s\in S$.
\end{proof}

Next, we try to prove that a set with exactly two elements is a fixed point of $f$.
Suppose that $ \min([\![k, k+1]\!]^f)=a$ and $\max([\![k, k+1]\!]^f)=b$.
 The following lemma shows that the maximum   and the minimum   of $X^f$ have a  close relation with   $k, a, b$ for  $X\in     \mathcal{P}_{\rm fin}(S)$.
\begin{lemma}  For any $X\in     \mathcal{P}_{\rm fin}(S)$, we have $$ \left\{\begin{array}{ccc}\min(X^f)&=&(a-k) (\max X-\min X)+\min X,\\
 \max(X^f)&=&(b-k)(\max X-\min X)+ \min X.\end{array}\right. $$
\end{lemma}
\begin{proof} If $X=\{r\}$ has a single element, since $ \{r\} ^f =  \{r\}$, the assertion is obvious .

Hereafter, we suppose $X$ has at least two elements.
Let $l$ be a positive integer. It follows from $[\![kl,kl+l]\!]=l[\![k, k+1]\!]$ that $$[\![kl, kl+l]\!]^f =l[\![k, k+1]\!]^f,$$
which leads to \begin{eqnarray}  \left\{\begin{array}{ccc}\min([\![kl, kl+l]\!]^f)&=& a l,\\
\max([\![kl, kl+l]\!]^f)&=&  b l.\end{array}\right.\end{eqnarray}
Since $(l-1) \{k\}^f+[\![k,k+l]\!]^f=[\![kl,kl+l]\!]^f$, by comparing the minimum and the maximum   of both sides we have
 \begin{eqnarray}  \left\{\begin{array}{ccc}\min([\![k, k+l]\!]^f)&=& (a-k)l+ k,\\
\max([\![k,k+l]\!]^f)&=&  (b-k)l+k.\end{array}\right.\end{eqnarray}
For $i\geq k$ and $l\geq 1$,  it follows from $$k [\![i,i+l]\!]^f=[\![ki,ki+kl]\!]^f=(i-1)\{k\}^f+[\![k,k+kl]\!]^f$$ that
\begin{eqnarray}  \left\{\begin{array}{ccc}\min([\![i, i+l]\!]^f)&=& (a-k)l+i,\\
\max([\![i, i+l]\!]^f)&=& (b-k)l+i.\end{array}\right.\end{eqnarray}
Furthermore, for a member $Y\in \mathcal{P}_{\rm fin}(S)$ with $\min Y=i\geq k$ and $\max Y=i+l$, since
 $Y+[\![i,i+l]\!]=[\![2i,2i+2l]\!]$, we have \begin{eqnarray}  \left\{\begin{array}{ccc} \min(Y^f)&=&\min([\![2i,2i+2l]\!]^f)-\min([\![i,i+l]\!]^f),\\ \max(Y^f)&=&\max([\![2i,2i+2l]\!]^f)-\max([\![i,i+l]\!]^f).\end{array}\right.\end{eqnarray}
 Applying equality (5), we have \begin{eqnarray}  \left\{\begin{array}{ccc} \min(Y^f)=(a-k)l+i ,\\ \max(Y^f)=(b-k)l+i .\end{array}\right.\end{eqnarray}

Now, suppose $X$ is an arbitrary member in $ \mathcal{P}_{\rm fin}(S)$ with at least two elements. Note that  $\min(X^f)=\min((X+\{k\})^f)-k$ and   $\max(X^f)=\max((X+\{k\})^f)-k$. Since $\min(X+\{k\})\geq k$, it follows from equality (7) that
 \begin{eqnarray}  \left\{\begin{array}{ccc} \min(X^f)=(a-k)(\max X-\min X)+\min X,\\ \max (X^f)=(b-k)(\max X-\min X)+ \min X.\end{array}\right.\end{eqnarray}
 \end{proof}

Applying Lemma 2.2, one can establish a relation between each member in $\{a,b\}$ and $k$.
\begin{lemma} Let   $ \min([\![k, k+1]\!]^f)=a,$ $\max([\![k, k+1]\!]^f)=b$.
Then $a=b-1=k$.
\end{lemma}
\begin{proof} Firstly, we prove $k\leq a$ and $k\leq b$. Suppose for a contradiction that $k=a+q$ with $q\geq 1$. Let $l$ be a positive integer with $l>\frac{k}{q}$.
Then it follows from Lemma 2.2 that $\min([\![k,k+l]\!]^f)=(a-k)l+k<0$, which is a contradiction and thus $k\leq a$.
Similarly, we  have $k\leq b$.

Secondly, we prove $a=k$. Suppose for a contradiction that $a>k$. For any given $  Y\not=\{\alpha_S\}\in \mathcal{P}_{\rm fin}(S)$  with $\min Y=\alpha_S$,  we have $$\min(Y^f)=(a-k)(\max Y-\min Y)+ \alpha_S>\alpha_S.$$
Since $f$ is a bijective mapping on $\mathcal{P}_{\rm fin}(S)$, there exists some $X$ with $\min X>\alpha_S$ such that $X^f=\{\alpha_S, k+1\}$. Then it follows from equality (8) that
$$\alpha_S=\min(X^f)=(a-k)(\max X-\min X)+\min X\geq\min X>\alpha_S,$$
which is a contradiction and thus $a=k$.

Finally, we  prove $b=k+1$. With $a=k$ in hands, equality (8) can be rewritten as  \begin{eqnarray}  \left\{\begin{array}{ccc} \min(X^f)&   =& \min X, \ \ \ \ \ \ \ \ \ \ \ \ \ \ \ \ \ \ \ \ \ \ \ \ \ \ \ \   \ \ \ \ \   \\   \max(X^f)&=&(b-k) (\max X-\min X)+ \min X,\end{array}\right.\end{eqnarray} for any $X\in \mathcal{P}_{\rm fin}(S)$. Suppose $Z_0 \in \mathcal{P}_{\rm fin}(S)$ satisfies $Z_0^f=[\![k,k+1]\!]$.
Then it follows from (9) that  \begin{eqnarray}  \left\{\begin{array}{ccccc} k&=&\min(Z_0^f)&   =&\min(Z_0),  \    \ \ \ \ \ \ \ \ \ \ \ \ \ \ \ \ \ \ \ \ \ \ \ \ \ \ \ \ \ \ \ \ \\  k+1&=& \max(Z_0^f)&=&(b-k) (\max(Z_0)-\min(Z_0))+ \min(Z_0),\end{array}\right.\end{eqnarray}
Equality (10) implies that  \begin{eqnarray}(b-k) (\max(Z_0)-\min(Z_0))=1.\end{eqnarray}
 Thus $b-k=\max(Z_0)-\min(Z_0)=1$ and $b=k+1$. \end{proof}

The application of Lemma 2.2 and Lemma 2.3 immediately gives the following corollary.
\begin{coro} Let $f$ be an automorphism of $\mathcal{P}_{\rm fin}(S)$ and let $k=\theta_S$.
Then \\
(i)\  \ $\min(X^f)=\min X$ and   $\max(X^f)=\max X$ for any $X\in \mathcal{P}_{\rm fin}(S)$;\\
(ii) \ $[\![i,j]\!]^f=[\![i,j]\!]$ for any positive integers $i,j$ with $k \leq i<j$.
\end{coro}
\begin{proof} (i) is a straightforward consequence of Lemma 2.2 and Lemma 2.3.

 Proof of  (ii):\ \ Following from (i), we have $\min([\![k,k+1]\!]^f)=k$ and $\max([\![k,k+1]\!]^f)=k+1$, which implies that $[\![k,k+1]\!]^f=[\![k,k+1]\!].$
Also, $\{q\}^f=\{q\}$ for any $q\in S$.

 For any positive integer $l$, since $ l [\![k,k+1]\!]  =  \{k(l-1)\} +[\![k,k+l]\!],$ both $[\![k,k+1]\!]$ and $\{k(l-1)\}$ are respectively fixed by $f$,
 by applying $f$ on two sides of the above equality, we have $  l [\![k,k+1]\!] = \{k(l-1)\}+[\![k,k+l]\!]^f$, which leads to $[\![k,k+l]\!]^f=[\![k,k+l]\!]$.

For any positive integers $i,j$ with $k \leq i\leq j$, since  $\{k\}^f= \{k\}, \{i\}^f=\{i\}$  and $[\![k,k+j-i]\!]^f=[\![k,k+j-i]\!]$,  by applying $f$ on two sides of $ \{k\}+[\![i,j]\!]  = \{i\}+  [\![k,k+j-i]\!]$, we have $$\{k\}+[\![i,j]\!]^f  = \{i\}+ [\![k,k+j-i]\!]=  [\![k+i,k+j]\!],$$ which leads to $[\![i,j]\!]^f=[\![i,j]\!].$ The proof is completed.
\end{proof}

Given $X \subset \mathbb{Z}$,   denote by $\Delta(X)$ the {\sf gap set} of $X$, i.e., the set of all integers $d \geq 1$
such that $\{x, x + d\} = X \cap [\![x, x + d[\!]$ for some $x \in \mathbb{Z}$. Accordingly, we define $g(X):={\rm sup}\ \Delta(X)$, called the   the {\sf maximal gap} of $X$ (see \cite{tri2}). If $X=\emptyset$, then we  put $ g(X)=0$.
For example, the maximal gap of $X=\{0,3,5,6,8,13\}$  is $5$.

\begin{lemma}  Let $f$ be an automorphism of $\mathcal{P}_{\rm fin}(S)$. Then $g(X^f)=g(X)$ for any $X\in \mathcal{P}_{\rm fin}(S)$.
\end{lemma}
\begin{proof} Assume that $g(X)=l+1$, and let $k=\theta_S$. It is easy to see that $$X+[\![k,k+l]\!]=[\![k+\min X,k+l+\max X]\!].$$ Applying $f$ we have $$X^f+[\![k,k+l]\!]=[\![k+\min X,k+l+\max X]\!],$$
which implies that $g(X^f)\leq l+1=g(X).$ By considering $f^{-1}$ on $X^f$, we have $g(X)=g((X^f)^{f^{-1}})\leq  g(X^f)$. Hence, $g(X^f)=g(X)$.
\end{proof}
An application of Lemma 2.5  tells us that a set $X\in \mathcal{P}_{\rm fin}(S)$ with  two elements is fixed by an automorphism of $\mathcal{P}_{\rm fin}(S)$.
\begin{lemma}  Let $f$ be an automorphism of $\mathcal{P}_{\rm fin}(S)$. If  $X\in \mathcal{P}_{\rm fin}(S)$ contains exactly two elements, then $X^f=X$.
\end{lemma}
\begin{proof} Suppose $X=\{x_1,x_2\}\in \mathcal{P}_{\rm fin}(S)$ with $x_1<x_2$. By Lemma 2.4, we have $\min(X^f)=x_1, \max(X^f)=x_2$. Since $g(X^f)=g(X)=x_2-x_1$ (by Lemma 2.5), $X^f$ has exactly two elements. Consequently, $X^f=\{x_1,x_2\}=X$.
\end{proof}

Let $S_\theta$ be the subsemigroup of $S$ consisting of all $l\in S$ for which $l\geq \theta_S$. That is to say that $S_\theta$ is the discrete interval  $[\![\theta_S, \infty)$.
One can easily see that the restriction of an automorphism of $\mathcal{P}_{\rm fin}(S)$ to $\mathcal{P}_{\rm fin}(S_\theta)$ is an automorphism of $\mathcal{P}_{\rm fin}(S_\theta)$.
\begin{lemma}  Let $f$ be an automorphism of $\mathcal{P}_{\rm fin}(S)$.  Then the restriction of $f$ to $\mathcal{P}_{\rm fin}(S_\theta)$ is an automorphism of $\mathcal{P}_{\rm fin}(S_\theta)$.
\end{lemma}
\begin{proof} Let $X\in S_\theta$. Since $\min(X^f)=\min X$ (by Corollary 2.4), we know $X^f\in  S_\theta$, thus the restriction of $f$ to $\mathcal{P}_{\rm fin}(S_\theta)$  is well-defined.
Clearly, the restriction of $f$ is an automorphism of $\mathcal{P}_{\rm fin}(S_\theta)$.
\end{proof}
For an automorphism $f$ of $\mathcal{P}_{\rm fin}(S)$, when we concentrate on what form does $f$ have, a key step is to investigate what form does the restriction of $f$ to $\mathcal{P}_{\rm fin}(S_\theta)$ have.
This reminds us to study (in the following section) the automorphisms of $\mathcal{P}_{\rm fin}(S_\theta)$.

\section{Automorphisms of the power semigroup of a discrete interval}
\quad\quad Following   the above section, when we   characterize the form of an automorphism $f$  on $\mathcal{P}_{\rm fin}(S)$,  a key step  is to investigate how $f$ acts on $\mathcal{P}_{\rm fin}(S_\theta)$,
where $S_\theta=[\![\theta_S, \infty)$ is the subset of $\mathbb{N}$ of all integers larger than or equal to $\theta_S$.  This is the reason why we turn to the automorphism group of  $\mathcal{P}_{\rm fin}(S_\theta)$ in this section.
For simple, suppose $\theta_S=k$ and thus $S_\theta=[\![k, \infty)$. Note that $k$ possibly is $0$ and if this case happens, then  $S_\theta=S$ and $\mathcal{P}_{\rm fin}(S_\theta)$ is a monoid with $\{0\}$ as zero-element.

Let $f$ be an automorphism of $\mathcal{P}_{\rm fin}(S_\theta)$. All properties about automorphisms of $\mathcal{P}_{\rm fin}(S)$  obtained in the above section still apply for $f$ (view $S_\theta$ as $S$).
We gather them   for later use.
\begin{lemma} Let $S_\theta=[\![k, \infty)$   with $X$ a non-empty finite subset, and    $f$ be an automorphism of  $\mathcal{P}_{\rm fin}(S_\theta)$. Then \\
(i) \ $\min(X^f)=\min X$ and $\max(X^f)=\max X$;\\
(ii) \ $[\![i,j]\!]^f=[\![i,j]\!]$ for any positive integers $i,j$ with $k \leq i\leq j$;\\
(iii) \   $g(X^f)=g(X)$;\\
(iv)\  $X^f=X$ if $X$ has at most two elements.
\end{lemma}

We continue study how $f$ acts on $\mathcal{P}_{\rm fin}(S_\theta)$.

\begin{lemma} Let   $f$ be an automorphism of  $\mathcal{P}_{\rm fin}(S_\theta)$ and let $X_0=\{k,k+2,k+3\}$. Then, either $X_0^f=X_0$, or  $X_0^f= \{k,k+1,k+3\}$.
\end{lemma}
\begin{proof} Applying (i) of Lemma 3.1, we know $\min(X_0^f)=\min(X_0)=k$ and $\max(X_0^f)=\max(X_0)=k+3$.
By (iii) of Lemma 3.1 we have $g(X_0^f)=g(X_0)=2$, which implies that  $X_0^f=X_0$ or $X_0^f= \{k,k+1,k+3\}$.
\end{proof}

Lemma 3.2 tells us that an automorphism $f$ possibly maps $ X_0=\{k,k+2,k+3\}$ to  $\{k,k+1,k+3\}$, which exactly is $\max(X_0)-X_0+\min(X_0)$.
This guides  us to  define an automorphism of   $\mathcal{P}_{\rm fin}(S_\theta)$ in such way:
\begin{lemma} Let   $\sigma$ be a transformation  on $\mathcal{P}_{\rm fin}(S_\theta)$ defined as $X\mapsto \max X-X+\min X,\   {\rm for\ all}\  X\in \mathcal{P}_{\rm fin}(S_\theta)$. Then $\sigma$ is an automorphism of $\mathcal{P}_{\rm fin}(S_\theta)$.
\end{lemma}
\begin{proof}  Clearly, $\sigma$ is an involution  (self-inverse transformation) on  $\mathcal{P}_{\rm fin}(S_\theta)$.
Moreover, for any $X, Y\in \mathcal{P}_{\rm fin}(S_\theta)$, $$\begin{array}{ccc}(X+Y)^\sigma&=&\max(X+Y)-(X+Y)+\min(X+Y)\ \ \ \ \ \ \ \ \ \   \\
&=&(\max X-X+\min X)+(\max Y-Y+\min Y)\\
&=&X^\sigma+Y^\sigma, \ \ \ \ \ \ \ \ \ \ \ \ \ \ \ \ \ \ \ \ \ \ \ \ \ \ \ \ \ \ \ \ \ \ \ \ \ \ \ \ \ \ \end{array}$$
which proves that $\sigma$ is an automorphism of  $\mathcal{P}_{\rm fin}(S_\theta)$.
\end{proof}

We believe that $\sigma$ is the only non-trivial automorphism of $\mathcal{P}_{\rm fin}(S_\theta)$.
The technique for proving the assertion is to reduce the   automorphisms   of  $\mathcal{P}_{\rm fin}(S_\theta)$ to  those of  $\mathcal{P}_{\rm fin,0}(\mathbb{N})$.
Recall that the  automorphisms   of    $\mathcal{P}_{{\rm fin},0}(\mathbb{N})$  have been obtained by   Tringali and Yan (Theorem 1.1).
 In order to apply  Theorem 1.1,  a preliminary work of us is to establish a relation between an automorphism  $f$ of $\mathcal{P}_{\rm fin}(S_\theta)$ and an automorphism  of  $\mathcal{P}_{\rm fin,0}(\mathbb{N})$.
We begin with a binary relation $``\sim"$ on  $\mathcal{P}_{\rm fin}(S_\theta)$.
\begin{defi} For $X,Y\in   \mathcal{P}_{\rm fin}(S_\theta)$, we write $X\sim Y$  if there is an integer $m\in \mathbb{Z}$ such that   $Y=m+X$.
\end{defi}
Clearly,  $``\sim"$ is an equivalence  relation. For $X\in   \mathcal{P}_{\rm fin}(S_\theta)$, set $$\overline X=\{Y: Y \in  \mathcal{P}_{\rm fin}(S_\theta), Y\sim
X\},$$ which is the equivalency class of $X$.
Denote by $\overline{\mathcal{P}_{\rm fin}(S_\theta)}$ the set of all    $\overline X$ for $X\in  \mathcal{P}_{\rm fin}(S_\theta)$.
Define an additive operation on $\overline{\mathcal{P}_{\rm fin}(S_\theta)}$ in a natural way: $$\overline X+\overline Y=\overline{X+Y}, \  {\rm for\ all}\ X,Y\in  \mathcal{P}_{\rm fin}(S_\theta).$$
It is easy to see that $$\overline{ X+ Y}=\overline{X_1+Y_1}  \ \ {\rm if}\ \ \overline X= \overline{X_1}, \ \overline{Y}=\overline{Y_1},$$
which implies the additive operation on $\overline{\mathcal{P}_{\rm fin}(S_\theta)}$ is well-defined. Thus $(\overline{\mathcal{P}_{\rm fin}(S_\theta)}, +)$ turns out to be a monoid with zero-element $\overline{\{k\}}$.

Let us investigate how an automorphism  of $\mathcal{P}_{\rm fin}(S_\theta)$ induces an automorphism of  $\overline{\mathcal{P}_{\rm fin}(S_\theta)}$.
For an automorphism $f$ of $\mathcal{P}_{\rm fin}(S_\theta)$, define $\overline f$ on  $\overline{\mathcal{P}_{\rm fin}(S_\theta)}$
in the way:$$   \overline X ^{\overline f}= \overline{X^f},\ \ {\rm for\ all}\  X\in \mathcal{P}_{\rm fin}(S_\theta).$$
For the sake that  $\overline f$ is well defined, we need the following property about $f$.
\begin{lemma} Let   $f$ be  an automorphism of $\mathcal{P}_{\rm fin}(S_\theta)$.
Then $(m+X)^f=m+X^f$ for all $m\in \mathbb{N}$ and $X\in  \mathcal{P}_{\rm fin}(S_\theta)$.\end{lemma}
\begin{proof} Considering the action of $f$ on $k+m+X$, we have $$(k+m+X)^f=\{k+m\}^f+X^f, \ \ (k+m+X)^f=\{k\}^f+(m+X)^f $$
Since $\{k+m\}^f=\{k+m\}$ and $\{k\}^f=\{k\}$  (by (iv) of Lemma 3.1), we have $$ k+m +X^f= k+ (m+X)^f,$$
which leads to $(m+X)^f=m+X^f$.
\end{proof}
\begin{lemma} Let   $f$ be  an automorphism of $\mathcal{P}_{\rm fin}(S_\theta)$.
Then the mapping $\overline f:  \overline X\mapsto \overline{X^f},\ \ {\rm for \ all}\ X\in \mathcal{P}_{\rm fin}(S_\theta) $,  is an automorphism of  $\overline{\mathcal{P}_{\rm fin}(S_\theta)}$. \end{lemma}
\begin{proof} With Lemma 3.5 in hands, it is easy to see that $\overline f$ is well defined. Indeed, if $\overline X=\overline{X_1}$, then there is some $m\in \mathbb{N}$ such that $X_1=m+X$ (or $X=m+X_1$).
 By Lemma 3.5, $X_1^f=m+X^f$ (or $X^f=m+X_1^f$) and $\overline{X_1^f}=\overline{X^f}$. Consequently, $\overline{X_1}^{\overline{f}}=\overline{X}^{\overline{f}}$, proving that  $\overline f$ is well defined.

 For the injectivity,   if $\overline{X_1}^{\overline{f}}=\overline{X}^{\overline{f}}$, then  $\overline{X_1^f} =\overline{X^f}$, and thus $X_1^f=m+X^f$ (or $X^f=m+X_1^f$) for some $m\in \mathbb{N}$.
 By Lemma 3.5,  $X_1^f=(m+X)^f$ (or $X^f=(m+X_1)^f$). Since $f$ is injective, $X_1=m+X$ (or $X=m+X_1$), which leads to  $\overline{X_1} =\overline{X}$, proving that $\overline f$ is injective.
 Obviously, $\overline f$ is surjective.

 Furthermore, $$\begin{array}{ccccccc}  (\overline X+\overline Y)^{\ \overline f}& =& (\overline{X+  Y})^{\ \overline f}&= &  \overline{(X+Y)^f}\\
    &=&  \overline{X^f+Y^f}&=&\overline{X^f}+\overline{Y^f}&=&\overline{X}^{\ \overline{f}}+\overline{Y}^{\ \overline{f}},\end{array} $$
  which implies that $\overline f$ is an automorphism of   $\overline{\mathcal{P}_{\rm fin}(S_\theta)}$.
\end{proof}
In order to reduce an automorphism $f$ of $\mathcal{P}_{\rm fin}(S_\theta)$ to that of $\mathcal{P}_{\rm fin,0}(\mathbb{N})$, we need further establish a relation between  $\overline{\mathcal{P}_{\rm fin}(S_\theta)}$ and  $\mathcal{P}_{\rm fin,0}(\mathbb{N})$.
Let $\varphi$ be the mapping from $\overline{\mathcal{P}_{\rm fin}(S_\theta)}$ to $\mathcal{P}_{\rm fin,0}(\mathbb{N})$ defined by:
$$\varphi: \overline X\mapsto X-\min X,  \ {\rm for\  all}\  X\in \mathcal{P}_{\rm fin}(S_\theta).$$
\begin{lemma}
The mapping $\varphi$ just defined is an isomorphism from   $\overline{\mathcal{P}_{\rm fin}(S_\theta)}$ to  $\mathcal{P}_{\rm fin,0}(\mathbb{N})$. \end{lemma}
\begin{proof} If $\overline X=\overline Y$, assume  (without loss of generality) that $X=m+Y$ for some $m\in \mathbb{N}$, which leads to  $\min X=m+\min Y$.
Consequently, $$\begin{array}{ccccc} \overline X  ^{\ \varphi}&=&X-\min X&=&(m+Y)-(m+\min Y)\\
&=&Y-\min Y&=&\overline Y^{\ \varphi},\ \ \ \ \ \ \ \ \ \ \ \ \ \ \  \ \ \ \ \ \ \ \ \ \end{array}$$
 proving that $\varphi$ is well defined.

For a given  $Y\in   \mathcal{P}_{\rm fin,0}(\mathbb{N})$,  $\overline{k+Y}$ lies in $ \overline{\mathcal{P}_{\rm fin}(S_\theta)}$ whose image under $\varphi$ is $Y$. This implies that $\varphi$ is surjective.

If $\overline{X_1}^{\ \varphi}=\overline{X}^{\ \varphi}$, then  $X_1-\min(X_1)=X-\min X$, and thus $\overline{X_1}=\overline{X}$, proving   that $\varphi$ is injective.

Furthermore, $$\begin{array}{ccc}  (\overline X+\overline{X_1})^{\ \varphi}& =&  \overline{X+ X_1} ^{\ \varphi}=   (X+X_1)-\min(X+X_1)\ \ \ \ \ \ \ \ \\
    &=&  (X-\min X)+(X_1-\min(X_1))=\overline{X}^{\ \varphi}+\overline{X_1}^{\ \varphi},\end{array} $$
  which implies that $\varphi$ is an isomorphism from   $\overline{\mathcal{P}_{\rm fin}(S_\theta)}$ to  $\mathcal{P}_{\rm fin,0}(\mathbb{N})$.
\end{proof}
 Based on Theorem 1.1, the automorphisms of ${\mathcal{P}_{\rm fin}(S_\theta)}$ are characterized as follows.
 \begin{theorem} Let  $S_\theta= [\![k,\infty)$  with $k\in \mathbb{N}$.  Then the only   automorphism of $\mathcal{P}_{\rm fin}(S_\theta)$ is either the identity or the
 the involution $\sigma: X \mapsto \max X -X+\min X\ {\rm for\ all}\  X\in \mathcal{P}_{\rm fin}(S_\theta)$.
\end{theorem}
 \begin{proof} Let $f$ be an automorphism of $\mathcal{P}_{\rm fin}(S_\theta)$ and let $\pi$ be the natural mapping from $\mathcal{P}_{\rm fin}(S_\theta)$ to $\overline{\mathcal{P}_{\rm fin}(S_\theta)}$ defined as
 $X^\pi= \overline X  \ {\rm for\ all}\  X\in \mathcal{P}_{\rm fin}(S_\theta)$.
 $$\begin{array}{ccccccc}\mathcal{P}_{\rm fin}(S_\theta)&\overrightarrow{\ \ \ \ \pi\ \ \ \ }&\overline{\mathcal{P}_{\rm fin}(S_\theta)}&\overrightarrow{\ \ \ \ \ \varphi\ \ \ \ }&\mathcal{P}_{\rm fin,0}(\mathbb{N})\\
 \downarrow  f& & \downarrow   \overline{f}& & \ \ \ \ \ \ \ \downarrow  \varphi^{-1}\overline f\varphi\\
 \mathcal{P}_{\rm fin}(S_\theta)&\overrightarrow{\ \ \ \ \pi\ \ \ \ \ }&\overline{\mathcal{P}_{\rm fin}(S_\theta)}&\overrightarrow{\ \ \ \ \ \varphi\ \ \ \ }&\mathcal{P}_{\rm fin,0}(\mathbb{N})\end{array}$$
 By the definition for $\pi, \overline f$ and $\varphi$, it is easy to check that the above diagram is commutative.

 $$\begin{array}{ccccccc}X&\overrightarrow{\ \ \ \ \pi\ \ \ \ }&\overline{X}&\overrightarrow{\ \ \ \ \ \varphi\ \ \ \ }&X-\min X\\
 \ \ \ \downarrow  f& & \ \ \ \downarrow   \overline{f}& & \ \ \ \downarrow  \varphi^{-1}\overline f\varphi\\
 X^f&\overrightarrow{\ \ \ \ \pi\ \ \ \ \ }&\overline{X^f}&\overrightarrow{\ \ \ \ \ \varphi\ \ \ \ }&X^f-\min X\end{array}$$
 For the automorphism $\overline f$ of $\overline{\mathcal{P}_{\rm fin}(S_\theta)}$, there exists an automorphism, say $\rho$, of $\mathcal{P}_{\rm fin,0}(\mathbb{N})$ such that $\overline f=\varphi\cdot \rho\cdot \varphi^{-1}$.
 Theorem 1.1 tells us $\rho$ is either the identity or the involution $Y\mapsto \max Y-Y$ for all $Y\in \mathcal{P}_{\rm fin,0}(\mathbb{N})$.
 If $\rho$ is the identity automorphism, then  $\overline f=\varphi\cdot \rho\cdot \varphi^{-1}$ is the identity automorphism of  $\overline{\mathcal{P}_{\rm fin}(S_\theta)}$ such that $\overline{X^f}=\overline X$ for all $X\in\mathcal{P}_{\rm fin}(S_\theta)$. Thus, for each given $X\in\mathcal{P}_{\rm fin}(S_\theta)$, we have $X^f=X+m$ (or $X=X^f+m$) with $m$ a non-negative integer. As $\min(X^f)=\min X$, we have $m=0$. In this case, $f$ is the identity automorphism of $\in\mathcal{P}_{\rm fin}(S_\theta)$.
 If $\rho$ is the involution $\sigma_0:Y\mapsto \max Y-Y \ {\rm for\ all}\  Y\in  \mathcal{P}_{\rm fin,0}(\mathbb{N})$. Then $$\overline{X^f}=(\overline X)^{\overline f}=(\overline X)^{\varphi\cdot\sigma_0\cdot \varphi^{-1}}=\overline{(X-\min X)^{\sigma_0}+k}.$$
 Then there exists an integer $m$ such that  \begin{eqnarray} \begin{array}{ccc}X^f&=&(X-\min X)^{\sigma_0}+k+m\ \  \ \ \ \ \ \ \  \ \ \ \ \ \  \ \ \ \ \ \ \ \ \\ &=& \max(X-\min X)-(X-\min X)+k+m\\ &=&\max X-X+k+m.\ \ \ \ \ \ \ \ \ \ \ \ \ \ \ \ \ \ \ \ \ \ \ \ \ \ \ \end{array}\end{eqnarray}
 Since $\min(X^f)=\min X$, by comparing the minimum   in $X^f$ and that in $\max X-X+k+m$, we have $k+m=\min X$. Consequently, $X^f=\max X-X+\min X$ for all $X\in \mathcal{P}_{\rm fin}(S_\theta)$, as required.
  \end{proof}

 \section{Proof of Theorem 1.3}

 \quad \quad Let $S$ be a numerical semigroup with $k=\theta_S$ be its critical element, $f$ an automorphism of $\mathcal{P}_{\rm fin}(S)$.
  If $k$ is the minimum   of $S$, then $S=S_\theta=[\![k,\infty)$ is a discrete interval. In this case, $f$ is either the identity automorphism, or the involution automorphism $\sigma:X\mapsto \max X-X+\min X$ for any $X\in \mathcal{P}_{\rm fin}(S)$.

 For the case when $k\not=\alpha_S$, by Lemma 2.7, $f|_{\mathcal{P}_{\rm fin}(S_\theta)}$, the restriction of $f$ to $\mathcal{P}_{\rm fin}(S_\theta)$, is an automorphism of $\mathcal{P}_{\rm fin}(S_\theta)$, which is either the identity or the involution automorphism $\sigma$.
 We claim the later case fails to happen. Since $\alpha_S\in S$ and $k-1\notin S$, there is $m\in[\![\alpha_S, k-2]\!]$ such that $m\in S$ and $m+1\notin S$. If $f|_{\mathcal{P}_{\rm fin}(S_\theta)}=\sigma$, then $$\begin{array}{ccccc}\{m,k,k+1\}^f+\{k\} &=& \{m+k,2k,2k+1\}^f
 = \{m+k,2k,2k+1\}^{\sigma}\ \ \ \ \ \ \ \ \ \ \ \ \ \ \ \\&=&\{m+k,m+k+1,2k+1\} =\{m, m+1,k+1\} +\{k\},\end{array}$$
 which leads to $\{m,k,k+1\}^f=\{m,m+1,k+1\}$ and $m+1\in S$, a contradiction to the choice of $m$.
 Consequently,  $f|_{\mathcal{P}_{\rm fin}(S_\theta)}$ is the identity automorphism of ${\mathcal{P}_{\rm fin}(S_\theta)}$, which fixes each $Y\in {\mathcal{P}_{\rm fin}(S_\theta)}$.
 Now, for each $X\in \mathcal{P}_{\rm fin}(S)$, we have  $X+\{k\}=(X+\{k\})^f=X^f+\{k\}$, proving that $X^f=X $ for all $ X\in \mathcal{P}_{\rm fin}(S)$.
 Hence, $f$ is the identity automorphism of  $\mathcal{P}_{\rm fin}(S)$, which completes the proof for Theorem 1.3.

\vskip 2mm
\noindent
{\bf Acknowledgement}

  We admire the referee for his (or her) deep insight into the topic of  this article, and we are grateful to the referee for his (or her) detailed reading and helpful suggestions, which  improve the readability of this article a lot!

\end{document}